\documentclass[12pt,reqno]{amsart}

\usepackage[latin1]{inputenc}
\usepackage[T1]{fontenc}
\usepackage{amssymb,url}

\theoremstyle{plain}
\newtheorem{theorem}{Theorem}[section]
\newtheorem{lemma}[theorem]{Lemma}
\newtheorem{proposition}[theorem]{Proposition}

\theoremstyle{definition}

\theoremstyle{remark}
\newtheorem{remark}[theorem]{Remark}
\newtheorem*{acknowledgements}{Acknowledgements}

\numberwithin{equation}{section}

\newcommand{\abs}[1]{\lvert #1\rvert}
\newcommand{\inp}[3][]{\left\langle #2,#3\right\rangle_{#1}}
\newcommand{\norm}[2][]{\left\lVert #2\right\rVert_{#1}}
\DeclareMathOperator{\ad}{ad}
\DeclareMathOperator{\diag}{diag}
\DeclareMathOperator{\grad}{grad}
\DeclareMathOperator{\Hom}{Hom}
\newcommand{\hook}{\mathbin{\lrcorner}}
\DeclareMathOperator{\Id}{Id}
\DeclareMathOperator{\im}{Im}
\newcommand{\Overline}[1]{\overline{\mathchoice{\displaystyle #1}{\textstyle
#1}{\scriptstyle #1}{\scriptscriptstyle #1}}}  
\DeclareMathOperator{\proj}{\mathbb P}
\DeclareMathOperator{\rank}{rank}
\DeclareMathOperator{\spec}{Spec}
\DeclareMathOperator{\Tr}{Tr}

\newcommand{\ar}[1]{\mathop{\rightleftarrows}%
\limits^{\alpha_{#1}}_{\beta_{#1}}}
\newcommand{\hqq}{\mathord{/\mkern-6mu/\mkern-6mu/}}

\newcommand{\LieG}{G}
\newcommand{\Liegc}{G^{\mathbb C}}
\newcommand{\lieg}{\mathfrak g}
\newcommand{\liegc}{\mathfrak g^{\mathbb C}}

\newcommand{\g}{} 
\def\g#1(#2){\mathfrak{#1}(#2)}
\newcommand{\G}{} 
\def\G#1(#2){\mathsf{#1}(#2)}
\newcommand{\gc}{} 
\def\gc#1(#2){\mathfrak{#1}(#2,\mathbb C)}
\newcommand{\Gc}{} 
\def\Gc#1(#2){\mathsf{#1}(#2,\mathbb C)}

\newcommand{\Sqrt}[1]{{#1}^{1/2}}
\newcommand{\transp}[1]{{#1}^{\mathsf T}}
\newcommand{\Nilp}{{\mathcal O}}
\newcommand{\Nj}[1]{\Nilp_{(#1)}}

\newenvironment{spmatrix}{\left(\smallmatrix}{\endsmallmatrix\right)}
\newcommand{\lzero}{\text{\large 0}}

\begin{document}
\title[HyperKähler Potentials]{HyperKähler Potentials via\\ 
Finite-Dimensional Quotients}

\author{Piotr Kobak}
\address[Kobak]{Instytut Matematyki\\
Uniwersytet Jagiello\'nski\\
ul.\ Rey\-mon\-ta~4\\
30-059 Krakow\\
Poland}
\email{kobak@im.uj.edu.pl}

\author{Andrew Swann}
\address[Swann]{Department of Mathematics and Computer Science\\
University of Southern Denmark\\
Odense University\\
Campusvej 55\\
DK-5230 Odense M\\
Denmark} 
\email{swann@imada.sdu.dk}


\begin{abstract}
  It is known that nilpotent orbits in a complex simple Lie algebra admit
  hyperKähler metrics with a single function that is a global potential for
  each of the Kähler structures (a hyperKähler potential).  In an earlier
  paper the authors showed that nilpotent orbits in classical Lie algebras
  can be constructed as finite-dimensional hyperKähler quotient of a flat
  vector space.  This paper uses that quotient construction to compute
  hyperKähler potentials explicitly for orbits of elements with small
  Jordan blocks.  It is seen that the Kähler potentials of Biquard and
  Gauduchon for $\Gc SL(n)$-orbits of elements with $X^2=0$, are in fact
  hyperKähler potentials.
\end{abstract}

\maketitle

\section{Introduction}

Adjoint orbits in complex semi-simple Lie algebras are known to carry a
compatible hyperKähler metric invariant under the compact group action (see
\cite{Kron90b,Kronheimer:semi-simple,Kova96,Biquard:orbits}).  Nilpotent
orbits are particularly interesting as they admit a hyperKähler structure
which is closely related to twistor spaces and quaternion-Kähler
geometries~\cite{Swann91b} and which comes equipped with a hyperKähler
potential.  If one only asks for a Kähler potential compatible with the
hyperKähler structure, then several examples are known.  Hitchin
\cite{Hitchin:integrable} gave an expression for a global Kähler potential
for a hyperKähler structure on the regular semi-simple orbit of~$\gc sl(n)$
in terms of theta functions.  Biquard and Gauduchon~\cite{BiGa96}
determined a simple formula for the Kähler potential for the hyperKähler
metric on semi-simple orbits of \emph{symmetric type}.  These orbits come
in continuous families and by taking a limit Biquard and Gauduchon also
obtain Kähler potentials for certain nilpotent orbits.

In \cite{KoSw99w,KoSw99hp}, Kähler and hyperKähler potentials were obtained
for orbits of cohomogeneity one and two by considering the invariants
preserved by the compact group action.  The cohomogeneity of a complex
orbit $\Nilp\subset\liegc$ is defined as the codimension of the generic
orbits of the compact group~$G$ on~$\Nilp$.  As the cohomogeneity
increases, we move further away from homogeneous manifolds and the geometry
of the orbits becomes more complicated.

But there are other ways of rating the level of complexity of nilpotent
orbits.  In the case when each simple component of $\liegc$ is classical
(i.e., equals $\gc su(n)$, $\gc so(n)$, or $\gc sp(n)$) it can be shown
that nilpotent orbits in $\liegc$ arise as hyperKähler reductions of the
flat hyperKähler spaces $\mathbb H^N$ (see~\cite{KoSw93a}).  This gives a
more explicit description of the hyperKähler metric and the corresponding
potential, as the latter comes simply from the radial function~$r^2$
on~$\mathbb H^N$.  The space $\mathbb H^N$ in the construction arises from
a \emph{diagram} of unitary vector spaces; the longer the diagram, the more
complicated the geometry of the orbit.  But even orbits that arise from the
simplest diagrams (i.e., those of length~$2$) may have arbitrary high
cohomogeneity, which puts them beyond the scope of the ``low cohomogeneity
approach'' mentioned above.  In \cite{KoSw93}, we successfully applied this
technique to construct the hyperKähler potential for the regular nilpotent
orbit in~$\gc sl(3)$, which has cohomogeneity~$4$.  The aim of this paper
is to apply the same construction to calculate hyperKähler potentials for
nilpotent orbits with diagrams of length two or three.  This includes
classical orbits of cohomogeneity one or two and also all orbits obtainable
as limits of semi-simple orbits of symmetric type.  In particular, we are
able to prove (in the $\gc sl(n)$ case) that the Kähler potentials obtained
by Biquard and Gauduchon on nilpotent orbits are in fact hyperKähler
potentials.  This is not apparent from their work, particularly because we
found in \cite{KoSw99hp} that several of these orbits admit families of
invariant hyperKähler metrics with Kähler potentials.  We also determine
the potential for orbits in $\gc so(n)$ which have length three diagrams
and Jordan type $(3,2^{2k},1^\ell)$.  In the simplest cases there is a
striking resemblance to the formulæ we have for the cohomogeneity two case,
but for $k\geqslant 2$ matters complicate rapidly.

In the calculations we use finite covering maps between nilpotent orbits
and the Beauville bundle construction.  It is worth pointing out that these
techniques combined with knowledge of the invariants of the compact group
action can be used to find the potential in several other cases, for
example for nilpotent orbits in the exceptional Lie algebra
$\lieg_2^{\mathbb C}$ (see \cite{KoSw99mg2}).

Explicit knowledge of hyperKähler potentials is of interest in the study of
real nilpotent orbits, cf.~\cite{Brylinski:instantons}, and we expect to
pursue this in future work.

The paper is organised as follows.  Section~\ref{sec:background} recalls
the hyperKähler quotient construction of classical nilpotent orbits and
gives some general results on hyperKähler potentials.  In section~3 we
derive formulæ for the potential for orbits with diagrams of length~$2$ and
then, in section~4, apply the result to the low cohomogeneity case.
Finally, in section~\ref{s:length3} we work out the potential for the
simplest orbits with diagrams of length~$3$.

\begin{acknowledgements}
  We are grateful for financial support from the \textsc{Epsrc} of Great
  Britain and \textsc{Kbn} in Poland.
\end{acknowledgements}

\section{Background and General Results}
\label{sec:background}

We begin by reviewing the general theory of the relationship between
hyperKähler quotients, hyperKähler potentials and nilpotent orbits.

A Riemannian manifold~$(N,g)$ with complex structures $I$, $J$ and $K$
satisfying the quaternion identities $IJ=K=-JI$, etc., is \emph{hyperKähler}
if $g$ is Hermitian with respect to each of the complex structures and the
two-forms $\omega_I(X,Y):=g(X,IY)$, $\omega_J$ and $\omega_K$ are closed.
Such a manifold is thus symplectic in three different ways.  If one
distinguishes the complex structure~$I$, then $N$ becomes a Kähler manifold
with a holomorphic symplectic two-form $\omega_c := \omega_J+i\omega_K$.

An interesting general problem is to find hyperKähler structures compatible
with a given complex structure~$I$ and a holomorphic symplectic
form~$\omega_c$.  One natural source of such manifolds is adjoint orbits
$\Nilp$ of a complex semi-simple Lie group~$\Liegc$.  Such an orbit
inherits a complex structure~$I$ as a submanifold of the complex vector
space~$\liegc$.  The complex symplectic form on~$\Nilp$ is given at
$X\in\Nilp$ by
\begin{equation*}
  \omega^\Nilp_c([A,X],[B,X]) = \inp X{[A,B]}, 
\end{equation*}
where $\inp\cdot\cdot$ is the \emph{negative} of the Killing form
on~$\liegc$.  If $\LieG$ is a compact real form of~$\Liegc$, then
$\Nilp$~admits a $\LieG$-invariant hyperKähler structure compatible with
$I$ and $\omega^\Nilp_c$
\cite{Kron90b,Kronheimer:semi-simple,Kova96,Biquard:orbits}.

The Marsden-Weinstein quotient construction was adapted to hyperKähler
manifolds in~\cite{Hitchin-KLR:hK}.  Suppose a Lie group $H$ acts on a
hyperKähler manifold~$N$ preserving $g$, $I$, $J$ and~$K$.  Suppose also
that there exist symplectic moment maps $\mu_I$, $\mu_J$ and $\mu_K$ from
$N$ to~$\mathfrak h^*$ for the action of $H$ with respect to the symplectic
forms $\omega_I$, $\omega_J$ and $\omega_K$.  For $I$, this means that for
each $V\in\mathfrak h$, the function $\mu_I^V:=\inp{\mu_I}V$ satisfies
\begin{equation}
  \label{eq:symplectic-moment}
  d\mu_I^V = \xi_V \hook \omega_I,
\end{equation}
where $\xi_V$ is the vector field generated by the action of~$V$.  We then
define a \emph{hyperKähler moment map} by
\begin{equation*}
  \mu\colon N\to \mathfrak h^*\otimes \im \mathbb H,\qquad
  \mu = \mu_Ii+\mu_Jj+\mu_Kk.
\end{equation*}
The \emph{hyperKähler quotient} of $N$ by~$H$ is defined to be 
\begin{equation*}
  N\hqq H := \mu^{-1}(0)/H.
\end{equation*}
If $H$ acts freely on~$N$, then $N\hqq H$~is a hyperKähler manifold of
dimension $\dim N-4\dim H$.  Even if the action of~$H$ is not free, there
is a natural way to write $N\hqq H$ as a union of hyperKähler manifolds
\cite{Dancer-Swann:geometry-qK}.  We will often distinguish the complex
structure~$I$ and write $\mu=(\mu_{\mathbb C},\mu_{\mathbb R})$, where
$\mu_{\mathbb C}=\mu_J+i\mu_K$ and $\mu_{\mathbb R}=\mu_I$.  The map
$\mu_{\mathbb C}$~is then a complex symplectic moment map for the
(infinitesimal) action of~$H^{\mathbb C}$ on~$N$.

For nilpotent orbits in the classical Lie algebras, a $\LieG$-invariant
hyperKähler metric may be constructed by finite-dimensional hyperKähler
quotients~\cite{KoSw93a}.  The only other orbits for which such a
construction is known are the semi-simple orbits in~$\gc sl(n)$
\cite{Nakajima:ale} together with finite quotients of a couple of orbits in
exceptional algebras~\cite{Kobak-Swann:exceptional}.  Let us briefly recall
the construction for nilpotent orbits.

\subsection{Nilpotent Orbits for Special Linear Groups}
\label{s:An:construction}

Given a nilpotent element $A\in\gc sl(n)$ such that $A^{k-1}\not=0$ and
$A^k=0$ one defines the associated \emph{image flag} to be $
\{0\}=V_0\rightleftarrows V_1\rightleftarrows V_2
\rightleftarrows\dots\rightleftarrows V_k=\mathbb C^n $, where $V_i=\im
A^{k-i}$.  We consider the complex vector space
\begin{equation}
  \label{f:def:W}
  W=\bigoplus_{i=0}^{k-1}
  \bigl(\Hom(V_i,V_{i+1})\oplus \Hom(V_{i+1},V_i)\bigr)
\end{equation}
and represent elements $(\dots,\alpha_i,\beta_i,\dots)$ of~$W$ by diagrams
\begin{equation*}
  \{0\}=V_0\ar0 V_1\ar1 V_2\ar2\cdots\ar{k-1} V_k=\mathbb C^n.
\end{equation*}

Taking $\mathbb C^n$ to be equipped with a Hermitian two-form, induces
Hermitian inner products on each $V_i$, $i=0,1,2,\dots,k$, and we get a
norm on~$W$ given by
\begin{equation}
  \label{radial:W}
  r^2=\norm{(\dots,\alpha_i,\beta_i,\dots)}^2
  =\sum_{i=1}^{k-1}\Tr (\alpha_i^*\alpha_i+\beta_i\beta_i^*).
\end{equation}
The inner products enables us to make sense of Hermitian adjoints
$\alpha_i^*$ and $\beta_i^*$ and to endow the vector space~$W$ with a
quaternionic structure by defining $j(\dots,\alpha_i,\beta_i,\dots)=
(\dots,-\beta_i^*,\alpha_i^*,\dots)$.

The product $H=\G U(V_1)\times\dots\times \G U(V_{k-1})$ of unitary groups
acts in a natural way on~$W$:
\begin{multline*}
  (a_1,\dots,a_{k-1})(\dots,\alpha_i,\beta_i,
  \dots,\alpha_{k-1},\beta_{k-1})\\ 
  = (\dots,a_{i+1}\alpha_i a_i^{-1}, a_i \beta_i
  a_{i+1}^{-1},
  \dots, \alpha_{k-1} a_{k-1}^{-1}, a_{k-1} \beta_{k-1}). 
\end{multline*}
This action preserves the quaternionic structure on~$W$, and the hyperKähler
moment map $\mu=(\mu_{\mathbb C},\mu_{\mathbb R})$ is given by
\begin{equation}
  \label{e:mm}
  \begin{split}
    \mu_\mathbb C
    &=(\dots,\alpha_i \beta_i-\beta_{i+1} \alpha_{i+1},\dots),\\
    \mu_\mathbb R
    &=(\dots,\alpha_i \alpha_i^*-\beta_i^* \beta_i+
    \beta_{i+1} \beta_{i+1}^* -\alpha_{i+1}^*
    \alpha_{i+1},\dots).
  \end{split} 
\end{equation}
The hyperKähler quotient $W\hqq H$ is homeomorphic to the closure
$\overline\Nilp$ of the nilpotent orbit $\Nilp=\Gc SL(n) A$, which is a
singular algebraic variety.  The identification is induced by the map
$\psi\colon W\to\gc gl(n)$ given by
\begin{equation}
  \label{e:equiv-map}
  \psi(\dots,\alpha_i,\beta_i,\dots)=\alpha_{k-1}\beta_{k-1}.
\end{equation}
If $W_0\subset W$ denotes the open set where each $\alpha_i$ is injective
and each $\beta_i$ is surjective, then $\psi\colon W_0\hqq H\to\Nilp$ is a
diffeomorphism.  In fact, $\psi$~is the complex symplectic moment map for
the action of~$\Gc GL(n)$ on $W_0\hqq H$ and so the general theory of
moment maps implies that $\psi^*\omega^\Nilp_c$~agrees with the complex
symplectic structure on~$W_0\hqq H$.  Note that $j$ on~$W$ acts on $\Nilp$
by $\alpha_{k-1}\beta_{k-1}\mapsto -\beta_{k-1}^*\alpha_{k-1}^*$ which
agrees with the real structure $X\mapsto -X^*$ on~$\gc sl(n)$ defining the
Lie algebra of the compact group~$\G SU(n)$.

\subsection{Nilpotent Orbits in Orthogonal and Symplectic Algebras}
\label{s:BCDn:construction}
The above construction may be adapted to the remaining classical Lie
algebras $\gc so(n)$ and $\gc sp(n)$.  We start with a nilpotent
element~$A$ in the Lie algebra~$\liegc$ with $A^k=0$ and $A^{k-1}\ne0$.
Let $\delta$ be~$0$, if $\liegc=\gc so(n)$, or $1$, if $\liegc=\gc sp(n)$.
We consider the image flag
\begin{equation}
  \label{f:diag:BCD}
  \{0\}\rightleftarrows(V_1,\omega_1)\rightleftarrows(V_2,\omega_2)
  \rightleftarrows\dots\rightleftarrows(V_k,\omega_k)
  =(\mathbb C^n,\omega_k),
\end{equation}
where $\omega_i\colon V_i\times V_i\to \mathbb C$ are non-degenerate
bilinear forms satisfying
\begin{equation*}
  \omega_i(X,Y)=(-1)^{k-i+\delta}\omega_i(Y,X).
\end{equation*}
(This implies that $\dim V_i$ is even if $k-i+\delta$ is odd).  We denote
by $\cdot^\dagger$ the adjoint with respect to the forms~$\omega_i$ and
define Lie groups
\begin{equation*}
  H_i=\{A\in \G U(V_i):A^\dagger A=\Id_{V_i}\}.
\end{equation*}
Then $H_i$ is $\G Sp(V_i)$, if $k-i+\delta$ is odd, or $\G O(V_i)$, if
$k-i+\delta$ is even.

Take $H=H_1\times\dots\times H_{k-1}$ and let $W$ be the quaternionic
vector space as in formula~\eqref{f:def:W}.  The subspace $W^+\subset W$
defined by the equations
\begin{equation*}
  \beta_i=\alpha_i^\dagger, \qquad i=1,\dots, k-1,
\end{equation*}
is a quaternionic vector space.  The equations~\eqref{e:mm} define a
hyperKähler moment map for the action of~$H$ on~$W^+$.  Using the map
$\psi$ of~\eqref{e:equiv-map}, the hyperKähler quotient $W^+\hqq H$ may be
identified with the closure of the nilpotent orbit $H_k^{\mathbb C}
A\subset \mathfrak h_k^{\mathbb C}$.  Again, this identification is
compatible with the complex-symplectic form~$\omega^\Nilp_c$ and the real
structure. 

\subsection{HyperKähler Potentials}
\label{s:hkp}
A real-valued function $\rho\colon N\to\mathbb R$ on a hyperKähler
manifold~$N$ is called a \emph{hyperKähler potential} if $\rho$ is
simultaneously a Kähler potential for each of the Kähler structures
$(\omega_I,I)$, $(\omega_J,J)$ and $(\omega_K,K)$.  For $I$, this means
that $\omega_I=i\overline{\partial_I}\partial_I\rho$, or equivalently
\begin{equation*}
  \omega_I = -\tfrac12 dId\rho.
\end{equation*}
In general, $N$~will not admit a hyperKähler potential even locally.
Indeed, the existence of~$\rho$ implies that if we set
$\zeta=\tfrac12\grad\rho$ then $\{\zeta,I\zeta,J\zeta,K\zeta\}$ generates
an infinitesimal action of $\mathbb H^*\cong\mathbb R\times\G Sp(1)$ such
that
\begin{equation*}
  L_{I\zeta} g = 0,\quad L_{I\zeta}I = 0,\quad\text{and}\quad L_{I\zeta}J=2K,
\end{equation*}
with similar expressions for the action of $J\zeta$ and $K\zeta$, obtained
by permuting $(I,J,K)$ cyclically (see \cite{Swann91b,Boyer-GM:reduction}).

We need to know how hyperKähler potentials behave with respect to
hyperKähler quotients.  An indirect proof of a slightly weaker form of the
following result may be found in~\cite{Swann91b}.  Beware that the
hypotheses given in~\cite{Boyer-GM:reduction} are not quite strong enough.

\begin{theorem}
  \label{thm:potential-quotient}
  Let $(N,g,I,J,K)$ be a hyperKähler manifold admitting a hyperKähler
  potential~$\rho$.  Suppose a Lie group~$H$ acts freely and properly
  on~$N$ preserving $g$, $I$, $J$, $K$ and~$\rho$.  Suppose also that there
  is a hyperKähler moment map $\mu$ for the action of~$H$ on~$N$ and that
  $\mu$~is equivariant with respect to the infinitesimal action of~$\G
  Sp(1)$ defined by~$\rho$, meaning
  \begin{equation}
    \label{eq:equivariance}
    L_{I\zeta}\mu_I = 0, \quad L_{I\zeta}\mu_J = -2\mu_K, \quad \text{etc.}
  \end{equation}
  Then the function~$\rho$ induces a hyperKähler potential on the
  hyperKähler quotient $N\hqq H$.
\end{theorem}

\begin{proof}
  Let $i\colon \mu^{-1}(0)\hookrightarrow N$ be the inclusion and write
  $\pi\colon \mu^{-1}(0) \to Q:= N\hqq H$ for the projection.  The
  hyperKähler structure on the quotient is defined by the relations
  $\pi^*\omega_I^Q = i^*\omega_I$, etc.  In particular, at each $x\in
  \mu^{-1}(0)$ the tangent space to the fibre is spanned by the vector
  fields~$\xi_V$, for $V\in\mathfrak h$ and
  $\left(T_x\mu^{-1}(0)\right)^\bot=\{I\xi_V,J\xi_V,K\xi_V:V\in\mathfrak
  h\}$.  Thus if $Y\in T_x\mu^{-1}(0)$ is orthogonal to each $\xi_V$, then
  $IY$, $JY$ and $KY$ lie in $T_x\mu^{-1}(0)$ too.
  
  As $\rho$~is invariant under the action of~$H$, it descends to define a
  function $\rho_Q\colon Q\to\mathbb R$ satisfying $\pi^*\rho_Q =
  i^*\rho$.  This implies $\pi^* d\rho_Q = i^* d\rho$.  Now $d\rho$ is
  metric dual to~$2\zeta$, so $\zeta$~commutes with the action of~$H$, and
  we claim that $\zeta$ is tangent to $\mu^{-1}(0)$.
  
  The equivariance condition~\eqref{eq:equivariance} gives,
  \begin{equation*}
    2\mu_K^V = -L_{I\zeta}\mu_J^V = -I\zeta\hook(\xi_V\hook \omega_J)
    = \omega_K(\xi_V,\zeta),
  \end{equation*}
  using the $J$ version of~\eqref{eq:symplectic-moment}.  But now
  \begin{equation*}
    L_\zeta\mu_K^V = \zeta\hook d\mu_K^V = \zeta\hook(\xi_V\hook \omega_K)
    = \omega_K(\xi_V,\zeta) = 2\mu_K^V. 
  \end{equation*}
  Thus $L_\zeta\mu=2\mu$ and $\zeta$~preserves $\mu^{-1}(0)$.

  For $V\in\mathfrak h$, we have 
  \begin{equation*}
    g(\zeta,\xi_V) = \tfrac12d\rho(\xi_V) = \tfrac12 L_{\xi_V}\rho = 0,
  \end{equation*}
  as $\rho$~is $H$-invariant.  So $I\zeta$ is also tangent to
  $\mu^{-1}(0)$.  In particular, $i^*Id\rho=Ii^*d\rho$ and we have
  \begin{equation*}
    \pi^*(-\tfrac12 dId\rho_Q) = i^*(-\tfrac12 dId\rho) = i^*\omega_I =
    \pi^*\omega_I^Q,
  \end{equation*}
  so $\rho_Q$ is a Kähler potential for $\omega_I^Q$.  Similar
  computations apply for $J$ and $K$ and we have that $\rho_Q$ is a
  hyperKähler potential on $Q=N\hqq H$.
\end{proof}

For the flat hyperKähler spaces $W$ and $W^+$ introduced above, the
hyperKähler potential is given by the function~$r^2$ of
equation~\eqref{radial:W}.  A hyperKähler potential on $\Nilp=W_0\hqq
H\subset \gc sl(n)$ or $\Nilp=W_0^+\hqq H\subset \gc so(n)$ or $\gc
sp(n)$ is then given by the restriction of~$r^2$ to the zero set of the
hyperKähler moment map.

One can now ask whether this hyperKähler potential is any sense unique.  In
fact, one can answer such a question for nilpotent orbits in general.  The
following is an extension of an argument in~\cite{Brylinski:instantons}.

\begin{proposition}
  \label{prop:unique-hKp}
  Let $\LieG$~be a compact semi-simple Lie group and let $\sigma$~be the
  corresponding real structure on~$\liegc$.  Let $\Nilp \subset \liegc$ be
  a nilpotent orbit with the Kirillov-Kostant-Souriau complex symplectic
  structure $(I,\omega^\Nilp_c)$.  Suppose $(g,I,J,K)$ is a hyperKähler
  structure on~$\Nilp$ such that
  \textup(a\textup)~$\omega_J+i\omega_K=\omega^\Nilp_c$,
  \textup(b\textup)~$g$ is invariant under the compact group~$G$ and
  \textup(c\textup)~the structure admits a hyperKähler potential such that
  for the induced $\mathbb H^*$-action $j\in \mathbb H^*$~acts
  as~$\sigma|_\Nilp$.  Then the hyperKähler structure is unique.
\end{proposition}

\begin{proof}
  By averaging with the $G$-action we may assume that there is a
  $G$-invariant hyperKähler potential~$\rho$ on~$\Nilp$.  Let
  $\zeta=\tfrac12\grad \rho$, as above.  Then $L_\zeta\omega_I=2\omega_I$
  and $L_\zeta\omega^\Nilp_c=2\omega^\Nilp_c$, so
  \begin{equation*}
    \omega^\Nilp_c = \tfrac12 d(\zeta\hook \omega^\Nilp_c).
  \end{equation*}
  Note that as $\omega^\Nilp_c$ is a $(2,0)$-form,
  $\zeta\hook\omega^\Nilp_c$~is of type~$(1,0)$.
  
  However, as $\Nilp$~is nilpotent, the form~$\omega^\Nilp_c$ is exact in
  Dolbeault cohomology: $\omega^\Nilp_c=d\theta$, with $\theta_X([X,A])=\inp
  XA$, which is holomorphic and $\Liegc$-invariant.  Therefore
  $\theta-\tfrac12\zeta\hook\omega^\Nilp_c$~is closed.  But
  $H^1(\Nilp,\mathbb C)=0$, as for nilpotent orbits have finite fundamental
  groups.  So $\theta-\tfrac 12\zeta\hook\omega^\Nilp_c = df$, for some
  function $f\colon\Nilp\to \mathbb C$.
  
  Now $df$~is of type~$(1,0)$ and holomorphic.  It is also
  $\LieG$-invariant, as $\zeta$ commutes with~$G$.  Therefore we may
  average $f$ over the action of~$G$ to get a $G$-invariant holomorphic
  function $\tilde f$ satisfying $d\tilde f = \theta - \tfrac12
  \zeta\hook\omega^\Nilp_c$.  However, such a function is
  $\Liegc$-invariant and $\Liegc$~acts transitively on~$\Nilp$, so $\tilde
  f$ is constant and $\zeta\hook\omega^\Nilp_c=2\theta$.  Therefore, the
  $(1,0)$-part of~$\zeta$ agrees with the $(1,0)$ part of the Euler vector
  field on~$\Nilp$.  As both these vector fields preserve~$I$, we have that
  $\zeta$~equals the Euler vector field.
  
  We now have that the quotient of~$\Nilp$ by the $\mathbb C^*$-action
  generated by $\zeta$ and $I\zeta$ is the projectivised
  orbit~$\proj(\Nilp)$ with $\theta$~as its complex-contact structure and
  with real structure~$\sigma$.  By~\cite{Swann:HTwNil}, $\proj(\Nilp)$~is
  the twistor space of a unique quaternion-Kähler manifold~$M$ of positive
  scalar curvature and $\Nilp$~is the associated hyperKähler
  manifold~$\operatorname{\mathcal U}(M)$.  Thus the hyperKähler structure
  is uniquely determined.
\end{proof}

\section{Nilpotent Orbits with Diagrams of Length Two} 
\label{s:length2}

Assume that $\liegc$ is a classical complex simple Lie algebra and
$\Nilp\subset\liegc$ is an orbit of a rank $k$ nilpotent matrix
$X\in\Nilp\subset\gc sl(n)$ which satisfies $X^2=0$.  Then $X$ has Jordan
type $(2^k, 1^{n-2k})$.  Such orbits are precisely those that arise from
diagrams of length two:
\begin{equation*}
  \{0\} \rightleftarrows
  \mathbb C^k \mathop{\rightleftarrows}\limits^\alpha_\beta
  \mathbb C^n.
\end{equation*}
It follows from \S\ref{s:An:construction} that there exist
$\alpha\colon\mathbb C^2\to\mathbb C^n$ and $\beta\colon\mathbb C^n\to
\mathbb C^2$, such that $X=\alpha\beta$, with
\begin{equation}
  \label{Un:length2mmeqs}
  \beta \alpha=0\quad \text{and}\quad \beta \beta^*=\alpha^* \alpha.
\end{equation}
When $\mathfrak g=\g su(n)$ this is the full set of equations for~$\Nilp$.
If $\mathfrak g$ is either $\g o(n)$ or $\g sp(n)$, then we have
additionally
\begin{equation}
  \label{dagger}
  \beta=\alpha^\dagger.
\end{equation}
In all cases $\rank\alpha=\rank\beta=\rank X = k$, so $\alpha$ is injective
and $\beta$ is surjective.

We shall use the above equations to calculate the hyperKähler potential
$\rho$ on~$\Nilp$.  From Theorem~\ref{thm:potential-quotient} we know that
$\rho$ is the restriction of the radial function~$r^2$.
By~\eqref{radial:W} we have
\begin{equation}
  \label{potential:alphabeta}
  \rho = \Tr(\alpha^*\alpha+\beta\beta^*)
  = 2\Tr\alpha^*\alpha
  = 2\Tr\Lambda,
\end{equation}
where $\Lambda=\alpha^*\alpha=\beta\beta^*$.  Since $\Lambda$ is
self-adjoint, there exists an orthonormal basis $\{e_1,\dots,e_k\}$ for
$\mathbb C^k$ in which $\Lambda$ is diagonal,
\begin{equation*}
  \Lambda=\diag(\lambda_1,\lambda_2,\dots,\lambda_k).
\end{equation*}
Thus $\rho=2(\lambda_1+\dots+\lambda_k)$.

Note that
\begin{equation*}
  \inp{\beta^* e_i}{\beta^* e_j} = \inp{\beta\beta^* e_i}{e_j}=
  \inp{\Lambda e_i}{e_j} = \lambda_i\delta_{ij}.
\end{equation*}
In particular, ${\norm {\beta^*e_i}}^2=\lambda_i$.  But $\beta^*$
is injective, so $\lambda_i>0$ and $\{\beta^* e_1,\dots,\beta^* e_k\}$ is
an orthogonal basis for $\im\beta^*$.

Now consider the matrix~$X^*X$.  On $\im\beta^*$, we have $X^*X=
\Lambda^2$, since
\begin{equation*}
  X^*X\beta^* e_i
  = \beta^*\alpha^*\alpha\beta\beta^* e_i
  = \beta^*\Lambda^2 e_i
  = {\lambda_i}^2\beta^* e_i.
\end{equation*}
On the other hand, $(\im\beta^*)^\perp=\ker\beta$ and $X=\alpha\beta$, so
$X^*X$ vanishes on $(\im\beta^*)^\perp$.  As a result $X^*X$ has
eigenvalues ${\lambda_1}^2,\dots,{\lambda_k}^2$.  Writing $\spec
X^*X=\{\mu_1,\dots,\mu_r\}$ with $\mu_i$ distinct and of multiplicity $k_i$
we get

\begin{theorem}
  \label{pot:length2}
  Let $\Nilp$ be the adjoint orbit of a non-zero nilpotent matrix~$X$ in a
  complex classical Lie algebra, and assume that $X^2=0$. Then the
  hyperKähler potential for the canonical hyperKähler metric on~$\Nilp$ is
  given by the formula
  \begin{equation}
    \label{eq:pot:length2}
    \rho(X)=2\sum_{\mu_i\in\spec(X^*X)}k_i{\mu_i}^{1/2}
  \end{equation}
\end{theorem}

\begin{remark}
  \label{X*X_eigenvals:rem}
  The above formula can be obtained from~\eqref{potential:alphabeta} by
  explicitly solving \eqref{Un:length2mmeqs} and~\eqref{dagger} for a given
  nilpotent element~$X$.  For example consider orbits in $\gc sl(n)$.  Then
  $X$ is $\G U(n)$-conjugate to
  \begin{equation}
    \label{block-nilp-matr}
    M=
    \begin{pmatrix}
      0&A\\
      0&0\\
    \end{pmatrix}
    ,
  \end{equation}
  where $A=\diag(a_1,\dots,a_k)$ with $a_i$ real and positive.  To see this
  note that $X^*X$ determines a set of orthonormal eigenvectors
  $e_1,\dots,e_k$ with positive eigenvalues $\mu_1,\dots,\mu_k$.  Moreover,
  $\inp{Xe_i}{Xe_j}=\mu_i\delta_{ij}$, so ${\mu_i}^{-1/2}Xe_i$,
  $i=1,\dots,k$ are also orthonormal.  Since $X^2=0$ it follows that
  \begin{equation*}
    0=\inp{X^2e_i}{Xe_j}=\inp{Xe_i}{X^*Xe_j}
    =\mu_j\inp{Xe_i}{e_j}.
  \end{equation*} 
  In effect the vectors
  \begin{equation*}
    e_1,\dots e_k, 
    {\mu_1}^{-1/2}Xe_1,\dots,{\mu_k}^{-1/2}Xe_k
  \end{equation*} 
  form an orthonormal set.  Complete this to an orthonormal basis in
  $\mathbb C^n$.  In this basis $X$ has the required form, with
  $a_i=\Sqrt{\mu_i}$.
  
  It follows that $X$ is $\G SU(n)$-conjugate to $\lambda M$ for some
  $\lambda$ satisfying $\lambda\overline\lambda=1$.  The moment map
  equations~\eqref{Un:length2mmeqs} are now solved by
  \begin{equation*}
    \alpha=\lambda
    \begin{pmatrix}
      A^{1/2}\\0
    \end{pmatrix} 
    \quad
    \text{and}
    \quad
    \beta=\overline\lambda
    \begin{pmatrix}
      0& A^{1/2}
    \end{pmatrix}
    ,
  \end{equation*}
  where $A^{1/2}=\diag(a_1^{1/2},\dots a_k^{1/2})$.  In particular
  $A=\alpha^*\alpha=\beta\beta^*$.  We have $\spec(XX^*)=\spec(A^2)
  =\{{a_1}^2,\dots,{a_k}^2\}$ and, by \eqref{eq:pot:length2}
  \begin{equation*}
    \rho(X)=2\sum_{i=1}^k\abs{a_i}.
  \end{equation*}
  This agrees with the formula obtained in~\cite{BiGa96}.  There Biquard
  \&{} Gauduchon showed that this formula gives a Kähler potential for a
  hyperKähler structure on the nilpotent orbit.  This was done by
  considering the orbit in $\gc sl(n)$ as a limit of semi-simple orbits.
  However, we have now shown that the Biquard-Gauduchon Kähler potential is
  in fact a hyperKähler potential.
\end{remark}

\section{HyperKähler Potentials for Low Cohomogeneity Orbits}

In the simplest case $\Nilp$ is a minimal nilpotent orbit in a
classical Lie algebra. Such orbit arises from a length two diagram. Its
Jordan type is given in Table~\ref{tab:k1}.  Minimal orbits are
cohomogeneity one so any two elements $X, X'\in \Nilp$ are conjugate if
and only if $\norm X=\norm {X'}$.  It follows that for all $X\in\Nilp$ the
matrix $X^*X$ has only one non-zero eigenvalue, say $\lambda$, with
multiplicity $\kappa$.  Then, by \eqref{eq:pot:length2}
$\rho=2\kappa\Sqrt\lambda$, so $\rho^2=4\kappa^2\lambda$.  But $\Tr
X^*X=\kappa\lambda$, so
\begin{equation}
  \label{min-orb-pot}
  \rho^2=4\kappa\Tr X^*X, \qquad\text{where}\quad
  \kappa=
  \begin{cases}
    1& \text{for $\gc sl(n)$, $\gc sp(n)$,}\\
    2& \text{for $\gc so(n)$.}
    \end{cases}
\end{equation}

\begin{table}[tbp]
  \begin{center}
    \newcommand{\tablestrut}{\vrule height 14pt depth 6pt width 0pt}
    \begin{tabular}[t]{|c||c|c|c|}
      \hline
      \tablestrut Type&$\gc sl(n)$&$\gc so(n)$&$\gc sp(n)$\\
      \hline \hline
      \tablestrut
      Cohomogeneity~$1$&$(21^{n-2})$&$(2^21^{n-4})$&$(21^{2n-2})$\\
      \hline
      \tablestrut Cohomogeneity~$2$ &$(2^21^{n-4})$&$(31^{n-3})$,
      $(2^41^{n-8})$
      &$(2^21^{2n-4})$\\
      \hline
    \end{tabular}
    \medskip
    \caption{Nilpotent orbits of low cohomogeneity in classical Lie
    algebras}
    \label{tab:k1}
  \end{center}
\end{table}

One finds the multiplicity $\kappa$ simply by calculating $X^*X$ where $X$
is the block matrix $\begin{spmatrix}A&0\\0&0 \end{spmatrix}$ with
$A=\begin{spmatrix}0&0\\1&0 \end{spmatrix}$ for $\gc sl(n)$ and $\gc
sp(n)$, and $A=\begin{spmatrix}0&1&0\\-1&0&i\\0&-i&0 \end{spmatrix}$ for
$\gc so(n)$.

In fact the potential on a minimal nilpotent orbit in any complex simple
Lie algebra is equal to $\norm X=\sqrt{\Tr X^*X}$, up to a constant
multiplier, see for example \cite{KoSw99w}.

It is known that, with one exception, the next-to-minimal orbits in complex
semi-simple Lie algebras are precisely the cohomogeneity-two
orbits~\cite{DaSw97}.  The exception is the next-to-minimal nilpotent orbit
in~$\gc sl(3)$ which has cohomogeneity~$4$.  This case was dealt with in
\cite{KoSw93} while in \cite{KoSw99hp} hyperKähler potentials for
cohomogeneity-two nilpotent orbits were calculated: the latter were
expressed in terms of two invariants $\eta_1(X):=-K(X,\sigma X)$ and
$\eta_2(X):=\eta_1([X,\sigma X])$, where $K$ denotes the Killing form. In
our situation it will be more convenient to use the following two
invariants (which in fact are multiples of $\eta_1$ and $\eta_2$):
\begin{align*}
  c_1(X)&=\Tr XX^*,\\
  c_2(X)&=\Tr YY^*,\qquad \text{where $Y=[X,X^*]$.}
\end{align*}

\begin{theorem}
  \label{t:hkp:coh2}
  Let $\Nilp$ be a cohomogeneity-two nilpotent orbit in a classical Lie
  algebra. Then the hyperKähler potential for $\Nilp$ is given by the formula
  \begin{equation*}
    \rho^2=4\kappa c_1+4\kappa\sqrt{2{c_1}^2-\kappa c_2}
  \end{equation*}
  where $\kappa=1$ for $\gc sl(n)$ and $\gc sp(n)$, and $\kappa=2$ for $\gc
  so(n)$.
\end{theorem}

In the proof we shall consider the three classes of orbits which have
length two diagrams, and postpone the length three case
to~\S\ref{s:length3}.

\begin{proof}
  We use the notation of Remark~\ref{X*X_eigenvals:rem}.  Since $\Nilp$ is
  a cohomogeneity-two orbit, $X^*X$ has at most two different eigenvalues.
  By considering a matrix defined in \eqref{block-nilp-matr}, with $a_1,
  a_2$ arbitrary, and $a_3=\dots=a_k=0$ one finds that for a generic
  element $X$ in nilpotent orbits $\Nj{2,1^{n-k}}\subset\gc sl(n)$, and
  $\Nj{2,1^{2n-k}}\subset\gc sp(n)$ we have $\spec(X^*X)=\{\mu_1,\mu_2\}$
  where the eigenvalues $\mu_1$, $\mu_2$ have multiplicities $\kappa=
  k_1=k_2=1$.  An element $X$ of $\Nj{2^4,1^{n-8}}\subset\gc so(n)$ has, by
  Lemma~\ref{so:double-multipl} below, eigenvalues with even
  multiplicities.  But $X^*X$ has rank 4 so again
  $\spec(X^*X)=\{\mu_1,\mu_2\}$, this time with multiplicities $\kappa
  =k_1=k_2=2$. This can be verified by a direct calculation: a typical
  matrix in this orbit is conjugate to the matrix obtained by taking $X$ as
  in \eqref{block-nilp-matr} with $a_1=-a_k, a_2=-a_{k-1}$ arbitrary, and
  $a_3=\dots=a_{k-2}=0$; note that this is possible if we take the
  quadratic form which defines $\gc so(n)$ to be
  $\frac12(x_1x_n+x_2x_{n-1}+\dots+x_nx_1)$, cf.~\S\ref{ort2}.

  From~\eqref{eq:pot:length2} we have $ \rho=
  2\kappa({\mu_1}^{1/2}+{\mu_2}^{1/2})$.  The invariants $c_i$ are not
  difficult to compute in terms of $\mu_1$ and~$\mu_2$:
  \begin{equation*}
    c_1=\Tr XX^*=\kappa({\mu_1}+{\mu_2}),
  \end{equation*} 
  and, since $X^2=0$, we have
  \begin{equation*}
    \begin{split}
      c_2 &=\Tr ([X,X^*][X,X^*]^*)=\Tr(XX^*-X^*X)^2
      =2\Tr(X^*X)^2\\
      &=2\kappa({\mu_1}^2+{\mu_2}^2).
    \end{split}
  \end{equation*}
  Thus
  \begin{gather*}
    \rho=2\kappa({\mu_1}^{1/2}+{\mu_2}^{1/2}),\\ 
    c_1=\kappa(\mu_1+\mu_2)\quad\text{and}\quad
    c_2=2\kappa({\mu_1}^2+{\mu_2}^2)
  \end{gather*}
  which leads to the required formula for length two orbits.
  
  There is only one cohomogeneity 2 orbit with diagram of length greater
  than two, for proof in this case see~\S\ref{ss:3}.
\end{proof}

The above proof used the following lemma:

\begin{lemma}
  \label{so:double-multipl}
  If $X\in\gc so(n)$ then the non-zero eigenvalues for $X^*X$ have even
  multiplicities.
\end{lemma}

\begin{proof}
  We consider $\mathbb C^n$ with the standard quadratic and Hermitian
  forms, so that $\gc so(n)$ consists of skew-symmetric matrices, and
  $X^*=\transp {\overline X}$.  Let $J$ denote the $\mathbb R$-linear
  automorphism of $\mathbb C^n$, defined by the formula
  \begin{equation*}
    Jv=X^*\overline v.
  \end{equation*}

  Suppose $\lambda$ is a non-zero eigenvalue of~$X^*X$ and that $v$ is a
  corresponding eigenvector.  Now $\transp X=-X$, so $X^*=-\overline X$, and
  we get
  \begin{equation*}
    \begin{aligned}
      X^*XJv&=X^*XX^*\overline v   
      =X^*\Overline{X^* X v}
      &=\lambda X^*\overline v=
      \lambda Jv,
    \end{aligned}
  \end{equation*}
  since the eigenvalues of $X^*X$ are real.  Thus $Jv$ is also a
  $\lambda$-eigenvector of~$X^*X$.
  
  Note that $J^2v=X^*\Overline{X^*\overline v}=-X^*Xv=-\lambda v$.  It
  follows that $v$ and $Jv$ are linearly independent.  We conclude that
  $\lambda$-eigenvectors with $\lambda\not=0$ come in pairs $v$, $Jv$ which
  span $J$-invariant two-dimensional $\lambda$-eigenspaces.
\end{proof}

\section{Orbits with Diagrams of Length Three}
\label{s:length3}

The hyperKähler potential calculations for orbits that correspond to diagrams of
length three can be quite involved, and the result is known only in few
special cases.  One of the early results is the calculation of the hyperKähler 
potential for the generic orbit $\Nj{3}\subset\gc sl(3)$, given in
\cite{KoSw93}. The formula
\begin{equation*}
  \rho(X)=2\sqrt{(a^{2/3}+c^{2/3})^3+b^2}, \qquad \text{where}
  \qquad X=
  \begin{pmatrix}
    0&a&b\\0&0&c\\0&0&0
  \end{pmatrix} 
\end{equation*} 
was derived from moment map equation~\eqref{e:mm} for $\Nj3$. This seems to
be the most efficient formula; the attempts to write the potential for this
orbit in another language, for example in terms of Lie algebra invariants,
yield much more complicated results.  Note, however, that the regular orbit
in $\gc sl(3)$ is a three-to-one quotient of the minimal orbit in
$\mathfrak g_2^{\mathbb C}$, the potential in question is proportional to
the invariant $\sqrt{c_1}$ on $\mathfrak g_2^{\mathbb C}$.

In this section we shall consider nilpotent orbits in $\gc so(n)$ which
have a single Jordan block of size three.  For nilpotent orbits in $\gc
so(n)$ the Jordan blocks of even size come in pairs, so these orbits have
Jordan type $(3,2^{2k},1^{n-4k-3})$ and the corresponding diagram is
\begin{equation*}
  \{0\} \rightleftarrows \mathbb C \rightleftarrows \mathbb C^{2k+2}
  \rightleftarrows \mathbb C^n
\end{equation*}
We may assume that the orthogonal structures $\omega_1$ on $\mathbb C$ and
$\omega_3$ on $\mathbb C^n$, cf.~formula \eqref{f:diag:BCD}, are the
standard quadratic forms. In particular $\gc so(n)$ consists of
skew-symmetric matrices.

By \S\ref{s:BCDn:construction}, the orbit $\Nj{3,2^{2k}}\subset\gc so(n)$
is a hyperKähler quotient
\begin{equation*}
  \mathbb H^{(2k+2)(n+1)}\hqq (\Gc Sp(k)\times \mathbb Z_2)
\end{equation*}
and $\Nj{2^{2k+2}}\subset\gc so(n+1)$ is $\mathbb H^{(2k+2)(n+1)}\hqq \Gc
Sp(k)$. This indicates that there is a $\mathbb Z_2$-quotient map
$\Nj{2^{2k+2}}\to \Nj{3,2^{2k}}$.  Moreover, the hyperKähler potentials on
$\Nj{2^{2k+2}}$ and on $\Nj{3,2^{2k}}$ are restrictions of the radial
function $r^2$ on $\mathbb H^{(2k+2)(n+1)}$, so they are preserved by the
quotient map.

Now $\Nj{2^{2k+2}}$ is given by a diagram of length two, so one can use
Theorem \ref{pot:length2} to calculate the potential for $\Nj{2^{2k+2}}$,
and hence for $\Nj{3,2^{2k}}$. By making the inverse to the two-to-one
quotient map explicit one gets an algorithmic method of calculating the
hyperKähler potential on $\Nj{3,2^{2k}}$.  This is shown in the following
technical lemma.

\begin{lemma}
  \label{l:2:1map}
  Let $X\in\Nj{3,2^{2k}}$ and denote by $x\in\mathbb C^n$ the
  \textup(unique up to sign\textup) vector such that $X^2=x\transp x$.
  Then the hyperKähler potential $\rho$ on $\Nj{3,2^{2k}}$ is given by the
  formula
  \begin{equation*}
    \rho(X)=2\sum_{\mu_i\in\spec(X'{X'}^*)}k_i{\mu_i}^{1/2}
    \qquad\text{where}\qquad
    X'=
    \begin{pmatrix}
      X & x\\ -\transp x&0 
    \end{pmatrix}.
  \end{equation*}
\end{lemma}

\begin{proof}
  We begin by writing down the diagram for $\Nj{3,2^{2k}}$:
  \begin{equation*}
    \{0\}\rightleftarrows V_1\ar1V_2 \ar2V_3,
    \qquad \text{with}
    \quad V_1=\mathbb C, 
    \quad V_2=\mathbb C^{2k+2},
    \quad V_3=\mathbb C^n,
  \end{equation*}
  and the corresponding moment map equations 
  \begin{gather}
    \label{l3mmeqns1c}
    \beta_1\alpha_1=0,\\
    \label{l3mmeqns2c}
    \alpha_1\beta_1=\beta_2\alpha_2,\\
    \label{l3mmeqns1r}
    \beta_1{\beta_1}^*={\alpha_1}^*\alpha_1,\\
    \label{l3mmeqns2r}
    \alpha_1{\alpha_1}^*+\beta_2{\beta_2}^*=
    {\beta_1}^*\beta_1+{\alpha_2}^*\alpha_2.
  \end{gather}
  We also have
  \begin{equation*}
    \beta_i={\alpha_i}^\dagger,
  \quad\text{and}\quad
    X=\alpha_2\beta_2.
  \end{equation*}
  
  Consider now the diagram
  \begin{equation*}
    \{0\}\rightleftarrows V_2 \ar{}V_1\oplus V_3.
  \end{equation*}
  The moment map equations are
  \begin{gather}
    \label{l2mmeqnsc}
    \beta\alpha=0,\\
    \label{l2mmeqnsr}
    \beta\beta^*=\alpha^*\alpha
  \end{gather}
  and it is easy to see that the map
  \begin{equation*}
    (\alpha_1,\alpha_2) \mapsto \alpha=\alpha_1^\dagger\oplus\alpha_2
  \end{equation*}
  transforms the solutions of \eqref{l3mmeqns1c}--\eqref{l3mmeqns2r} into
  solutions of \eqref{l2mmeqnsc}--\eqref{l2mmeqnsr}. To verify this simply
  write $\alpha$ and $\beta$ in block-matrix form:
  \begin{equation*}
    \alpha=
    \begin{pmatrix}
      \alpha_2\\
      \beta_1
    \end{pmatrix}
    , \qquad
    \beta=\alpha^\dagger=
    \begin{pmatrix}
      \beta_2&-\alpha_1
    \end{pmatrix}.
  \end{equation*}
  Then it is clear that \eqref{l3mmeqns2c} is equivalent to
  \eqref{l2mmeqnsc} and \eqref{l3mmeqns1r} to \eqref{l2mmeqnsr}. The
  remaining two equations \eqref{l3mmeqns1c} and \eqref{l3mmeqns1r} are
  $\Gc O(1)=\mathbb Z_2$ moment map equations and are trivially satisfied.
  
  Note that if $(\alpha_1,\beta_1,\alpha_2,\beta_2)$ solves
  \eqref{l3mmeqns1c}--\eqref{l3mmeqns2r} then so does
  $(-\alpha_1,-\beta_1,\allowbreak \alpha_2,\allowbreak \beta_2)$.  This
  corresponds to
  \begin{equation*}
    \alpha=
    \begin{pmatrix}
      \alpha_2\\
      -\beta_1 \end{pmatrix}
    , \qquad
    \beta=\alpha^\dagger=
    \begin{pmatrix}
      \beta_2&\alpha_1
    \end{pmatrix}.
  \end{equation*}
  A solution $(\alpha_1,\beta_1,\alpha_2,\beta_2)$ represents an element
  $X=\beta_2\alpha_2\in \Nj{2^{2k+2}}$ while the lifts $X'_\pm$ are given
  by
  \begin{equation*}
    X'_\pm=\alpha\beta=
    \begin{pmatrix}
      \alpha_2\beta&\mp\alpha_2\alpha_1\\
      \pm\beta_1\beta_2&0
    \end{pmatrix}
  \end{equation*}
  Define $x=\alpha_2\alpha_1(1)$.  With our conventions ($\omega_1$ and
  $\omega_3$ are the identity matrices) the dagger operator acts on maps
  $\mathbb C\to\mathbb C^n$ as the transpose, so
  $\beta_1\beta_2=(\alpha_2\alpha_1)^\dagger=\transp x$. Also,
  \begin{equation*}
    \begin{split}
      X^2&=(\alpha_2\beta_2)^2=\alpha_2\beta_2\alpha_2\beta_2\\
      &=\alpha_2\alpha_1\beta_1\beta_2=x\transp x
    \end{split}
  \end{equation*}
  where the penultimate equality follows from \eqref{l3mmeqns2c}.  This
  shows that $X'$ is of the required form. Finally, note that
  \begin{equation*}
    \begin{split}
      {r'}^2&=\Tr \alpha\alpha^*+\Tr\beta^*\beta\\
      &=\Tr\alpha_2\alpha_2^*+\Tr\beta_2^*\beta_2+
      \Tr\alpha_1\alpha_1^*+\Tr\beta_1^*\beta_1=r^2
    \end{split}
  \end{equation*}
  which shows directly, that the two-to-one map respects the hyperKähler
  potentials.
\end{proof}

We shall apply the above lemma to determine the hyperKähler potential on
$\Nj{3,2^{2k}}$ in a few simple cases.  The first completes the proof of
Theorem~\ref{t:hkp:coh2}.

\subsection{$\Nj{3,1^{n-3}}$ in $\gc so(n)$.}
\label{ss:3}
As in Lemma~\ref{l:2:1map} define
\begin{equation*}
  X'= 
  \begin{pmatrix}
    X&x\\-\transp x&0
  \end{pmatrix}.
\end{equation*}
Then $X'$ lies in the minimal nilpotent orbit $\mathcal
N_{2^2,1^{n-3}}\subset\gc so(n+1)$, so the potential is given by
\eqref{min-orb-pot}.  We have
\begin{equation*}
  {X'}^*=
  \begin{pmatrix}
    X^*&-\overline x\\
    x^*&0
  \end{pmatrix},
\end{equation*}
and
\begin{equation}
  \label{l1+l2}
  \rho^2= 4\kappa\Tr X'{X'}^* = 4\kappa(\Tr XX^*+2\norm x^2).
\end{equation}
Putting $Y=[X,X^*]$ we get
\begin{equation*}
  \begin{split}
    c_2
    &:= \Tr YY^*\\
    &= \Tr(XX^*-X^*X)(XX^*-X^*X)\\
    &=2\Tr(XX^*)^2-2\Tr X^2{X^*}^2\\
    &=2\Tr(XX^*)^2-2\norm x^4
  \end{split}
\end{equation*}
since $X=x\transp x$.  

We know that $\rank X=2$ so $X^*X$ has at most two non-zero eigenvalues. It
follows from Lemma~\ref{so:double-multipl} that it has a unique non-zero
double eigenvalue, which we denote by $\lambda$.  Then, in a
suitable basis,
\begin{equation*}
  XX^*=\diag(\lambda,\lambda,0,\dots,0),
\end{equation*}
so $c_2=4\lambda^2-2\norm x^4={c_1}^2-2\norm x^4$, since $c_1=\Tr
XX^*=2\lambda$.  This implies that $\norm x^2=\sqrt{(c_1^2-c_2)/2}$.  Thus
\begin{equation*}
  \rho^2=4\kappa(c_1+2\norm x^2)=4\kappa c_1+4\kappa\sqrt{2c_1^2-2c_2}
\end{equation*}
which ends the proof of Theorem~\ref{t:hkp:coh2}.

\subsection{$\Nj{3,2^2,1^{n-7}}$ in $\gc so(n)$.}
\label{ort2}
For this orbit $X'{X'}^*$ has two double eigenvalues,
$\spec(X'{X'}^*)=\{\lambda_1,\lambda_2\}$, so the computation
of~\eqref{l1+l2} yields $\lambda_1+\lambda_2$ and not $\rho^2$:
\begin{equation}
  \label{newtr}
  2(\lambda_1+\lambda_2)=\Tr X'{X'}^*=\Tr XX^*+2\norm x^2=c_1+2\norm x^2
\end{equation}
Moreover, 
by Theorem \ref{pot:length2}, 
\begin{equation*}
  \rho^2=4({\lambda_1}+{\lambda_2} +2\sqrt{\lambda_1\lambda_2})
\end{equation*}
so one needs to calculate the product of eigenvalues. This can be done by
calculating $\Tr (X'{X'}^*)^2$ but then it is necessary to determine
invariants like $\norm{X \overline x}^2$. The most straightforward approach
is to take a generic nilpotent element $X$, augment it to get $X'$, and
find the eigenvalues of $X'$.

To simplify the calculations we can use the action of the compact group $\G
SO(n)$ on $\Nilp$ to put $X$ in a \emph{canonical form}. This is achieved
by using the Beauville bundle \cite{Beau97}. We shall briefly outline this
approach here; it is explained in more detail in
\cite[Section~4]{KoSw99hp}.

Consider $e\in\Nilp\subset\liegc$ and choose $f,h\in\liegc$ so that $e,f,h$
is an $\gc sl(2)$-triple. Then use the $\ad_h$-eigenspaces $\liegc(i)$ to
define the algebras
\begin{equation*}
  \mathfrak p=\bigoplus_{i\geqslant0}\liegc(i), \qquad\mathfrak
  n=\bigoplus_{i\geqslant2}\liegc(i).  
\end{equation*}
It turns out that $\mathfrak p$ is a parabolic algebra and it does not
depend on the choice of $f,h$. This gives what is sometimes referred to as
the \emph{canonical} fibration $\Nilp\to\mathcal F$ where $\mathcal
F=\Liegc/P$ is a flag manifold with $P$ the normaliser of $\mathfrak p$.
Moreover, $\Nilp$ is an open dense subset of the \emph{Beauville bundle}
\begin{equation*}
  N(\Nilp)=\Liegc\times_{\!P}\mathfrak n,
\end{equation*}
the canonical fibration being the restriction to $\Nilp$ of the Beauville
bundle fibration.

Choose a flag $v\in \mathcal F$. Since $\mathcal F$ is $G$-homogeneous any
element $e\in\Nilp$ can be moved by the action of the compact group $G$
into the Beauville bundle fibre $N(\Nilp)_v$. It is enough to calculate the
hyperKähler potential $\rho$ for nilpotent elements $e\in\Nilp\cap
N(\Nilp)_v$.

We now calculate the hyperKähler potential on $\Nj{3,2^2,1^{n-7}}$.  First,
assume that the quadratic form on $\mathbb C^n$ is given by the
anti-diagonal matrix $(S)_{ij}$ with $S_{ij}=\delta_{i,n+1-j}$.  Then $\gc
so(n)$ consists of matrices that are skew-symmetric about the
anti-diagonal.  The advantage of this choice for the quadratic form is that
nilpotent matrices in $\gc so(n)$ are $\Gc SO(n)$-conjugate to matrices
consisting of Jordan blocks.  In our situation we can arrange for the size
three block to be in the middle with the size two blocks placed
symmetrically about the anti-diagonal:
\begin{equation*}
  e=\begin{spmatrix}
    J_2&&0\\
    &J_3&\\
    0&&-J_2
  \end{spmatrix}
  \quad\text{with}\quad
  J_2=
  \begin{spmatrix}
    0&1\\
    0&0
  \end{spmatrix},
  \quad
  J_3=
  \begin{spmatrix}
    0&1&0\\
    0&0&-1\\
    0&0&0
  \end{spmatrix}.
\end{equation*}
(For simplicity we write everything for $\Nj{3,2^2}$, the formulæ
are identical in other cases.)

Moreover, we can choose the maximal torus to consist of the diagonal
matrices $\diag(a_1,a_2,\dots,-a_2,-a_1)$.  Then we have an $\gc
sl(2)$-triple $e,f,h$ with $e$ as above and $h=\diag(1,1,2,0-2,-1,-1)$. To
make matters simpler we use the Weyl group to rearrange the diagonal matrix
$h$, so take $h'=\diag(2,1,1,0,-1,-1,-2)$. It is enough to work out the
$\ad_{h'}$ eigenspaces that have eigenvalues $\geqslant 2$ to see that a
typical element of the Beauville bundle fibre has the following form
\begin{equation*}
  Y=
  \begin{spmatrix}
    &&&a&b&0&0\\    
    &&&0&v&0&0\\
    \smash{\lzero}&&&0&0&-v&-b\\
    &&&0&0&0&-a\\
    \lzero&&&&&\lzero
  \end{spmatrix}
  .
\end{equation*}
(To be precise the $(1,6)$ and $(6,1)$ entries in the matrix have weight 3
and thus belong to the Beauville bundle fibre, but one can assume they
vanish by using the action of the stabiliser $\G SO(2)\G SO(2)\G Sp(1)$.)

The aim is now to apply Lemma~\ref{l:2:1map} but we need to go back to the
standard basis, where the quadratic form is diagonal.  To diagonalise the
quadratic form $S$ consider the matrix $Q$ written in a block form:
\begin{equation*}
  Q=\tfrac1{\sqrt2}
  \begin{spmatrix}
    \mathbf 1_3&0&-i\mathbf 1_3\\
    0&\sqrt2&0\\
    \mathbf 1_3&0&i\mathbf 1_3
  \end{spmatrix}
  ,
\end{equation*}
where $\mathbf 1_3$ is the $3\times3$ identity matrix. Then $ \transp Q S Q
= 1, $ so $X=Q^{-1}YQ=Q^*YQ$ is skew-symmetric. Lemma \ref{l:2:1map},
applied to $X$, gives $x=\frac 1{\sqrt2}\transp{(ai,0,\dots,0,ai)}$, and a
direct calculation yields
\begin{equation*}
  \lambda_1\lambda_2=2\abs a^2\abs v^2.
\end{equation*}
Let us introduce a new invariant 
\begin{equation*}
  c_{21}=c_1(X^2)=\Tr XXX^*X^*=\norm {X^2}^2.
\end{equation*}
A simple calculation shows that
\begin{equation*}
  {c_1}^2-c_2-2c_{21}=8\abs a^2\abs v^2
  =4\lambda_1\lambda_2.
\end{equation*}
By combining this with \eqref{newtr} we get the following formula:
\begin{proposition}
  The hyperKähler potential $\rho$ for the canonical hyperKähler structure
  on the nilpotent orbit $\Nj{3,2^2}\subset \gc so(n)$ is given by the
  formula
  \begin{equation*}
    \rho^2=8c_1+16\sqrt{c_{21}}+16\sqrt{{c_1}^2-c_2-2c_{21}}.
  \end{equation*}
\end{proposition}

Note the similarity of this formula to that in Theorem \ref{t:hkp:coh2}:
for length two diagrams $c_{21}=0$ while in the cohomogeneity two situation
(the orbit $\Nj{3,1^{n-3}}$ in $\gc so(n))$ the invariant $c_{21}$ is a
combination of $c_1$ and $c_2$.

\subsection{$\Nj{3,2^4,1^{n-11}}$ in $\gc so(n)$.}

Finally, we shall only indicate here how the matters tend to complicate if
one tries to proceed in the same manner and calculate the hyperKähler
potential for $\Nj{3,2^4,1^{n-11}}$. We start with the same strategy as in
the previous section (again, it is enough to analyse the case of
$\Nj{3,2^4}$).

Here we take the semi-simple element
\begin{equation*}
  h'=\diag(2,1,1,1,1,0,-1,-1,-1,-1,-2).
\end{equation*}
Taking into account the action of the stabiliser $\G SO(2)\G SO(2)\G
Sp(2)$, a typical element of the fibre of the Beauville bundle can be
written as
\begin{equation*}
  Y=\begin{spmatrix}
    &&&a&b&0&0&0&0\\        
    &&&0& v_1&-w_2&w_3&0&0\\
    &&&0& v_2&w_1&0&-w_3&0\\
    \smash{\lzero}&&&0& v_3&0&-w_1&w_2&0\\
    &&&0&0&-v_3&-v_2&-v_1&-b\\
    &&&0&0&0&0&0&-a\\
    \lzero&&&&&&\lzero
  \end{spmatrix}
\end{equation*}
As before, set $X=Q^*YQ$ and then $x=\frac
1{\sqrt2}\transp{(ai,0,\dots,0,ai)}$. Calculations now become complex
enough and the authors used \textsc{Maple}.

The result can be described as follows.  Denote $v=\transp{(v_1,v_2,v_3)}$
and $w=\transp{(w_1,w_2,w_3)}$.  Also, write $\zeta=\transp v w=\sum
v_iw_i$.  Then the hyperKähler potential $\rho$ for the canonical
hyperKähler metric on $\Nj{3,2^4,1^{n-11}}$ is given by the formula
\begin{equation*}
  \rho=2({\lambda_1}^{1/2}+{\lambda_2}^{1/2}+{\lambda_3}^{1/2})
\end{equation*}
where $\lambda_i$ are the roots of the cubic $z^3-pz^2+qz-r$ with
\begin{align*}
  p&=2\abs a^2+\abs b^2+\abs v^2+\abs w^2=c_1+\abs a^2\\
  q&=\abs\zeta^2+\abs b^2\abs w^2+2\abs a^2(\abs v^2+\abs
  w^2)\\
  r&=\abs a^2\abs\zeta^2.
\end{align*}

\providecommand{\bysame}{\leavevmode\hbox to3em{\hrulefill}\thinspace}

\end{document}